\documentclass[12pt]{article}
\usepackage{amssymb,amsmath}

\begin{document}

\baselineskip15pt

\font\fett=cmmib10
\font\fetts=cmmib7
\font\bigbf=cmbx10 scaled \magstep2
\font\bigrm=cmr10 scaled \magstep2
\font\bbigrm=cmr10 scaled \magstep3
\def\cnl{\centerline}
\def\etal{{\it et al.}}
\def\text{{}}
\def\un{^{(n)}}
\def\ff{{\cal F}}
\def\Ref#1{(\ref{#1})}

\def\Blm{\left|}
\def\Brm{\right|}
\def\Bl{\left(}
\def\Br{\right)}
\def\nti{n\to\infty}
\def\lnti{\lim_{\nti}}
\def\law{{\cal L}}
\def\sji{\sum_{j\ge1}}
\def\Cal{\cal}

\def\real{\text{\rm I\kern-2pt R}}
\def\re{\real}
\def\expec{\text{\rm I\kern-2pt E}}
\def\ex{\expec}
\def\prob{\text{\rm I\kern-2pt P}}
\def\pr{\prob}
\def\rat{\text{\rm Q\kern-5.5pt\vrule height7pt depth-1pt\kern4.5pt}}
\def\comp{\text{\rm C\kern-4.7pt\vrule height7pt depth-1pt\kern4.5pt}}
\def\nat{\text{\rm I\kern-2pt N}}
\def\integ{{\bf Z}}
\def\qedbox{\vcenter{\hrule height.5pt\hbox{\vrule width.5pt height8pt
\kern8pt\vrule width.5pt}\hrule height.5pt}}
\def\half{{\textstyle {1 \over 2}}}
\def\quarter{{\textstyle {1 \over 4}}}
\def\scrl{{\Cal L}}
\def\scrd{{\Cal D}}
\def\tod{\buildrel \scrd \over \longrightarrow}
\def\eqd{\buildrel \scrd \over =}
\def\tolone{\buildrel L_1 \over \longrightarrow}
\def\toltwo{\buildrel L_2 \over \longrightarrow}
\def\tolp{\buildrel L_p \over \longrightarrow}
\def\l{\lambda}
\def\L{\Lambda}
\def\a{\alpha}
\def\b{\beta}
\def\g{\gamma}
\def\G{\Gamma}
\def\d{\delta}
\def\D{\Delta}
\def\e{\varepsilon}
\def\h{\eta}
\def\z{\zeta}
\def\th{\theta}
\def\k{\kappa}
\def\m{\mu}
\def\n{\nu}
\def\p{\pi}
\def\r{\rho}
\def\s{\sigma}
\def\S{\Sigma}
\def\t{\tau}
\def\f{\varphi}
\def\ch{\chi}
\def\ps{\psi}
\def\o{\omega}
\def\lee{\,\le\,}
\def\leee{\quad\le\quad}
\def\gee{\,\ge\,}
\def\geee{\quad\ge\quad}
\def\scrn{{\Cal N}}
\def\scra{{\Cal A}}
\def\scrf{{\Cal F}}
\def\var{\text{\rm Var\,}}
\def\cov{\text{\rm Cov\,}}

\def\sin{\sum_{i=1}^n}
\def\sn{\sum_{i=1}^n}
\def\cross{\times}

\def\nin{\noindent}
\def\bsk{\bigskip}
\def\msk{\medskip}
\def\widebar{\bar}

\def\sno{\sum_{n\ge 0}}
\def\sj{\sum_{j\ge 0}}
\def\proof{\noindent{\bf Proof.}\quad}
\def\remark #1. {\medbreak\noindent{\bf Remark #1.}\quad}
\def\remarks{\medbreak\noindent{\bf Remarks.}\quad}
\def\head #1/ {\noindent{\bf #1.}\medskip}
\def\dtv{d_{TV}}
\def\Po{\text{\rm Po\,}}
\def\Bi{\text{\rm Bi\,}}
\def\Be{\text{\rm Be\,}}
\def\CP{\text{\rm CP\,}}

\def\pb{\hbox{{\fett\char'031}}}
\def\pbh{\text{\bf{\hbox{\^{\fett\char'031}}}}}
\def\bp{\pb}
\def\lb{\hbox{{\fett\char'025}}}
\def\lbs{\hbox{{\fetts\char'025}}}
\def\bl{\lb}
\def\rb{\hbox{\fett\char'032}}
\def\sb{\hbox{\fett\char'033}}
\def\iid{\text{\rm independent\ and\ identically\ distributed}}

\def\and{\ \text{\rm{ and }}\ }
\def\for{\ \text{\rm{ for }}\ }
\def\forall{\ \text{\rm{ for all }}\ }
\def\forsome{\ \text{\rm{ for some }}\ }
\def\tif{\text{{\rm if\ }}}
\def\tin{\text{{\rm in\ }}}
\def\lcm{\text{\rm{l.c.m.\,}}}

\def\BCL{Barbour, Chen and Loh}
\def\ABT{Arratia, Barbour and Tavar\'e}
\def\BHJ{Barbour, Holst and Janson}

\def\Blb{\left\{}
\def\Brb{\right\}}

\def\giv{\,|\,}
\def\Giv{\,\Big|\,}
\def\ep{\hfill$\Box$}
\def\bone{{\bf 1}}
\def\non{\nonumber}

\def\adb{}
\def\sodmi{\sum_{s=0}^{d-1}}
\def\sie{\sum_{i\in E}}
\def\Le{\ \le\ }
\def\Eq{\ =\ }
\def\TP{{\rm TP\,}}
\def\Ord#1{O\Bl #1 \Br}
\newcommand{\eqs}{\begin{eqnarray*}}
\newcommand{\ens}{\end{eqnarray*}}

\newcommand{\eqa}{\begin{eqnarray}}
\newcommand{\ena}{\end{eqnarray}}
\newcommand{\eq}{\begin{equation}}
\newcommand{\en}{\end{equation}}

\def\numberlikeadb{\global\def\theequation{\thesection.\arabic{equation}}}
\numberlikeadb
\newtheorem{theorem}{Theorem}[section]
\newtheorem{lemma}[theorem]{Lemma}
\newtheorem{corollary}[theorem]{Corollary}
\newtheorem{proposition}[theorem]{Proposition}
\newtheorem{example}[theorem]{Example}

\def\pr{{\mathbb P}}
\def\ex{{\mathbb E}}
\def\re{{\mathbb R}}
\def\integ{{\mathbb Z}}
\def\ZZ{\integ}
\def\Ceka{\v Cekan\-avi\v cius}
\def\sjom{\sum_{j=0}^m}
\def\uo{^{(0)}}
\def\ul{^{(l)}}
\def\gg{{\mathcal G}}
\def\bh{{\bar h}}
\def\bX{{\widebar X}}
\def\bU{{\widebar U}}
\def\BX{Barbour \&~Xia}
\def\sro{\sum_{r\ge0}}
\def\srz{\sum_{r\in\integ}}
\def\ssz{\sum_{s\in\integ}}
\def\diff{\lfloor \m-\s^2 \rfloor}
\def\diffi{\lfloor \m_1-\s_1^2 \rfloor}
\def\difft{\lfloor \m_2-\s_2^2 \rfloor}
\def\sid{\s_1^2 + \d_1}
\def\std{\s_2^2 + \d_2}
\def\hX{{\widehat X}}
\def\hY{{\widehat Y}}
\def\hW{{\widehat W}}
\def\hQ{{\widehat Q}}
\def\hS{{\widehat S}}
\def\Ge{{\rm Ge\,}}
\def\sjn{\sum_{j=1}^n}
\def\tpnfg{[T_i^+-i > n/4]}
\def\tmnfg{[i-T_i^- > n/4]}
\def\tpnfl{[T_i^+-i \le n/4]}
\def\tmnfl{[i-T_i^- \le n/4]}
\def\bbar{{\bar\b}}
\def\tb{{\tilde \b}}

\title{Translated Poisson approximation\\ for Markov chains}
\author
{A.\ D.\ Barbour\thanks{Angewandte Mathematik,
Winterthurerstrasse~190,
CH--8057 Z\"URICH, Switzerland: {\tt a.d.barbour@math.unizh.ch}
} ~and Torgny
Lindvall\thanks{School of Mathematical Sciences, Chalmers and GU,
S--41296 G\"OTEBORG, Sweden: {\tt lindvall@math.chalmers.se}}
\\
Universit\"at Z\"urich and G\"oteborgs Universitet}
\date{}
\maketitle
\vglue-1cm
\begin{abstract}
The paper is concerned with approximating the distribution of a
sum~$W$ of integer valued random variables~$Y_i$, $1\le i\le n$,
whose distributions depend on the state of
an underlying Markov chain~$X$.  The approximation is in terms of
a translated Poisson distribution, with mean and variance chosen
to be close to those of~$W$, and the error is measured with
respect to the total variation norm.  Error bounds comparable to
those found for normal approximation with respect to the weaker
Kolmogorov distance are established, provided that the distribution
of the sum of the~$Y_i$'s between the successive visits of~$X$ to
a reference state is aperiodic.  Without this assumption,
approximation in total variation cannot be expected to be good.
\end{abstract}


\section{Introduction}\label{intro}
 \setcounter{equation}{0}
 The Stein-Chen method is now well
established in the study of approximation by a Poisson or compound
Poisson distribution (Arratia, Goldstein \& Gordon~(1990),
\BHJ~(1992)). It has turned out to be very efficient
for treating sums of the form $W := W_n := \sn Y_i$, where the
variables $Y_1, Y_2, \ldots$ are non-negative, integer-valued, rarely
different from 0, and have a short range of dependence.  A
basic example is the following: let $Y_1, Y_2, \ldots$ be independent
and taking values 0 or 1 only, with $p_i := {\mathbb P}(Y_i=1)$
generally small, to make
a Poisson approximation plausible. Then the method offers a proof of
the celebrated Le~Cam theorem, which is transparent and relatively
simple (\BHJ~1992, I.(1.23)), and gives the optimal constant:
\eq\label{LeCam}
\|\law(W) - \Po(\l)\| \le 2\l^{-1}\sn p_i^2
          \le 2\max_{1\le i\le n}p_i,
\en
where $\l := \ex W = \sn p_i$. Here, $\law(X)$ denotes the distribution
of a random element~$X$, $\Po(\l)$ the Poisson
distribution with mean~$\lambda$, and~$\|\n\|$ the
total variation norm of a signed
bounded measure~$\n$; we need this only for differences of probability
measures $Q,Q'$ on the integers~${\mathbb Z}$, when
\[
\|Q-Q'\|:=\sum_i |Q(i)-Q'(i)|  = 2 \sup_{A\subset\integ}|Q(A)-Q'(A)|.
\]

Clearly, if the~$p_i$'s are not required to be small, there is little
content in~\Ref{LeCam}. This is to be expected, since then $\ex W = \l$ and
$\var W = \l - \sn p_i^2$ need no longer be close to one another,
whereas Poisson distributions have equal mean and variance.
This makes it more natural to try to find a family of distributions for
the approximation within which both mean and
variance can be matched, as is possible using the normal family
in the classical central limit theorem.  One choice
is to approximate with a member of the family of
translated Poisson distributions $\{\TP(\m,\s^2),\,(\m,\s^2) \in
\re\times\re_+\}$, where
\eqs
\TP(\m,\s^2)\{j\} := \Po(\s^2 + \d)\{j - \lfloor \m-\s^2 \rfloor\}\\
  = \Po(\l')\{j - \g\},
  \quad j \in \integ,
\ens
where
\eqa
&&\g := \g(\m,\s^2) := \lfloor \m-\s^2 \rfloor,\quad
\d := \d(\m,\s^2) := \m-\s^2 - \g \non\\
&&\mbox{and}\quad
\l' := \l'(\m,\s^2) := \s^2 + \d.\label{dgl-def}
\ena
The $\TP(\m,\s^2)$ distribution is just that of a Poisson with
mean $\l' := \l'(\m,\s^2)$ \hbox{$:= \s^2 + \d$}, then shifted along
the lattice by an amount $\g := \g(\m,\s^2) := \lfloor \m-\s^2 \rfloor$.
In particular, it has mean $\l' + \g = \m$ and variance~$\l'$
\hbox{such that}  \hbox{$\s^2 \le \l' < \s^2+1$;}
note that $\l' = \s^2$ only if
\hbox{$\m-\s^2 \in \integ$}.  For sums of independent, integer-valued
random variables~$Y_i$, this idea has
been exploited by  Vaitkus \&~\Ceka~(1998), and also
in \BX~(1999), \Ceka\ \&~Vaitkus~(2001)
and Barbour \&~\Ceka~(2002), using Stein's method, leading to error
rates of the same order as in the
classical central limit theorem, but now with respect to the
much stronger total variation norm, as long as some `smoothness'
of the distribution of~$W$ can be established.

As in the Poisson case, the introduction of
Stein's method raises the possibility of making similar
approximations for sums of dependent random variables as well.
However, the `smoothness' needed is a bound of order
$O(1/\sqrt n)$ for $\|\law(W+1) - \law(W)\|$, entailing much
more delicate arguments than are required for Poisson approximation.
The elementary example of $2$--runs in independent Bernoulli
trials was treated in \BX~(1999), but the argument
used there was long and involved.
More recently, R\"ollin~(2005) has proposed
an approach which is effective in a wider range of
circumstances, including many local and combinatorial dependence
structures, in which one can find an imbedded sum of independent
Bernoulli random variables.
In this paper, we consider a different kind of dependence,
in which the distributions of the random variables~$Y_i$ depend on
an underlying Markovian environment.

We suppose that
$X=(X_i)^\infty_{i=0}$ is an
aperiodic, irreducible and stationary Markov chain with finite state
space $E=\{0,1,\ldots, K\}$. Let $Y_0, Y_1,\ldots$ be
integer-valued variables which are independent conditional on~$X$, and,
as in a hidden Markov model,
such that the conditional distribution $\law(Y_i\giv X)$ depends on
the value of~$X_i$ alone; we assume further that, for each $0 \le k\le K$,
the distributions $\law(Y_i \giv X_i=k)$ are the same for all~$i$.
Under these assumptions, and with $W = \sn Y_i$,
we show that $\|\law(W)-\TP(\ex W,\var W)\|$ is asymptotically small,
under reasonable conditions on the conditional distributions
$\law(Y_1 \giv X_1=k)$, $0\le k\le K$.
The detailed results are given in Theorems \ref{Th-1}--\ref{Th-3}.
Roughly speaking, we show that if these conditional distributions
are stochastically dominated by a distribution with
finite third moment, and if, as smoothness condition, the
distribution $Q := \law\Bl\sum_{i=1}^{S_1} Y_i \giv X_0=0\Br$
is aperiodic ($Q\{d\integ\} < 1$ for all $d\ge2$), where~$S_1$
is the step at which~$X$ first returns to~$0$, then
\eq\label{1.3}
\|\law(W)-\TP(\ex W,\var W)\| = \Ord{n^{-1/2}}.
\en
An ingredient of our argument, reflecting R\"ollin's~(2005)
approach, is again to find an appropriate imbedded sum of
independent random variables.

In the next section, we give an introduction to proving
translated Poisson approximation by way of the Stein--Chen method.
Lemma~\ref{Stein-approx} provides a generally applicable formula
for bounding the resulting error.  In Section~\ref{markov}, we
establish bounds on the total variation distance between $\law(W)$
and~$\law(W+1)$ using coupling arguments.  The results of these
two sections are combined in Section~\ref{main} to prove the main
theorems. Theorem~\ref{Th-3} gives rather general conditions
for~\Ref{1.3} to hold, whereas Theorem~\ref{Th-1}, in a
somewhat more restrictive setting, provides a relatively explicit formula
for the approximation error.
We then discuss the relationship of our results to those of
\Ceka\ \&~Mikalauskas~(1999), who studied the degenerate case in which
$Y_1 = h(k)$ a.s.\ on $\{X_1=k\}$, $0\le k\le K$.
\adb{We conclude by showing that, if~$Q$ is
in fact periodic, $\law(W)$ is usually not well approximated by a
translated Poisson distribution.}

\section{Translated Poisson approximation}
 \setcounter{equation}{0}
Since the $\TP(\m,\s^2)$ distributions are just translates of Poisson
distributions,  the Stein--Chen method can be used to establish total
variation approximation.  In particular, $W\sim \TP(\m,\s^2)$ if and only
if
\eq\label{A.1}
\ex\{\l'f(W+1) - (W-\g)f(W)\} = 0
\en
for all bounded functions $f:\integ\to\re$, where
$\l' = \l'(\m,\s^2)$ and  $\g = \g(\m,\s^2)$ are as defined
in~\Ref{dgl-def}.
Define $f^*_C$ for $C\subset\integ_+$ by
\[
\begin{split}
&f^*_C(k) = 0,\quad k\le0;\\
&\l'f^*_C(k+1) - kf^*_C(k) = \bone_C(k) - \Po(\l')\{C\},\quad k\ge0,
\end{split}
\]
as in the Stein--Chen method. It then follows that
\[
\|f^*_C\| \le (\l')^{-1/2}\ \mbox{ and }\ \|\D f^*_C\| \le (\l')^{-1}
\]
(\BHJ\ 1992, Lemma~I.1.1), where $\D f(j) := f(j+1)-f(j)$ and,
for bounded functions
$g:\integ\to\re$, we let $\|g\|$ denote the supremum norm.
Correspondingly, for $B\subset\integ$ such that
$B^* := B - \g \subset \integ_+$, the function~$f_B$ defined by
\eq\label{f_B-def}
f_B(j) := f^*_{B^*}(j-\g),\quad j\in\integ,
\en
satisfies
\eqa
\lefteqn{\l'f_B(w+1) - (w-\g)f_B(w)}\non\\
&=& \l'f^*_{B^*}(w-\g+1) - (w-\g)f^*_{B^*}(w-\g)\non\\
&=&\bone_{B^*}(w-\g) - \Po(\l')\{B^*\} \non\\
&=& \bone_B(w) - \TP(\m,\s^2)\{B\} \label{A.2}
\ena
if $w \ge \g$, and
\eq\label{A.3}
\l'f_B(w+1) - (w-\g)f_B(w) = 0
\en
if $w < \g$; and clearly
\eq\label{A.4}
\|f_B\| \le (\l')^{-1/2}\ \mbox{ and }\ \|\D f_B\| \le (\l')^{-1}.
\en
This can be exploited to prove the closeness in total variation of $\law(W)$
to $\TP(\m,\s^2)$ for an arbitrary integer-valued random variable~$W$.
The next two results make use of this.

\begin{lemma}\label{TP-diff}
Let $\m_1,\m_2\in\re$ and $\s_1^2,\s_2^2 \in \re_+\setminus\{0\}$ be such that
$\g_1 =\diffi \le \g_2 = \difft$. Then
\[
\|\TP(\m_1,\s^2_1) - \TP(\m_2,\s^2_2)\| \le 2\{\s_1^{-1}|\m_1-\m_2|
  + \s_1^{-2}(|\s_1^2-\s_2^2| + 1)\}.
\]
\end{lemma}

\proof
Both distributions assign probability~$1$ to $\integ \cap
\bigl[\g_1,\infty\bigr)$, so it suffices to consider~$B$ such that
$B - \g_1 \subset \integ_+$. Then, if $W \sim \TP(\m_2,\s_2^2)$, we have
\eqa
\lefteqn{\pr(W \in B) - \TP(\m_1,\s_1^2)\{B\}}\non\\
&=& \ex\{\bone_B(W) - \TP(\m_1,\s_1^2)\{B\}\}\non\\
&=& \ex\{\l_1f_B(W+1) - (W-\g_1)f_B(W)\},\non
\ena
from~\Ref{A.2}, where $\l_l := \l'(\m_l,\s_l^2)$, $l=1,2$.
Applying~\Ref{A.1}, it thus follows that
\eqs
\lefteqn{\pr(W \in B) - \TP(\m_1,\s_1^2)\{B\}}\\
&=& \ex\{(\l_1 - \l_2)f_B(W+1) - (\g_2 - \g_1)f_B(W)\}\\
&=& \ex\{(\l_1 - \l_2)\D f_B(W) - (\m_2-\m_1)f_B(W)\},
\ens
\vbox{\nin
and hence, from~\Ref{A.4}, that
\eqs
\lefteqn{|\pr(W \in B) - \TP(\m_1,\s_1^2)\{B\}|} \\
   &\le& (\l_1)^{-1}(|\s_1^2-\s_2^2| + |\d_1-\d_2|)
          + (\l_1)^{-1/2}|\m_1-\m_2|,
\ens
proving the lemma.
\ep
}

\bsk
The next lemma provides a very general means to establish total
variation bounds; it is our principal tool in Section~\ref{main}.
Note that we make no assumptions about the dependence structure
among the random variables $Y_1,\ldots,Y_n$.

\begin{lemma}\label{Stein-approx}
Let $Y_1,Y_2,\ldots,Y_n$ be integer valued random variables with finite
means, and define
$W := \sn Y_i$.  Let $(a_i)_{i=1}^n$ and $(b_i)_{i=1}^n$ be real
numbers such that, for all bounded $f:\integ\to\re$,
\begin{equation}\label{i-bd}
 |{\mathbb E}[Y_i f(W)]-{\mathbb E}[Y_i]{\mathbb E}f(W)-a_i
 {\mathbb E}[\Delta f(W)]|\leq b_i\|\Delta f\|, \quad 1\le i\le n.
\end{equation}
Then
\[
\|\law(W) - \TP(\ex W,\s^2)\| \le 2(\l')^{-1}\Bl \d + \sn b_i \Br
  + 2\pr[W < \ex W - \s^2],
\]
where $\s^2 := \sn a_i$, $\d = \d(\ex W,\s^2)$ and $\l' = \s^2 + \d$.
\end{lemma}

\proof
Adding~\Ref{i-bd} over~$i$, and then adding and subtracting
$c{\mathbb E}f(W)$ for $c \in \re$ to be chosen at will, we get
\[
|{\mathbb E}[(W-c)f(W)]-({\mathbb E}W-c-\s^2){\mathbb E}f(W)
 -\s^2 {\mathbb E}[f(W+1)]|\leq \Bl\sn b_i\Br\|\Delta f\|,
\]
where $\s^2=\sn a_i$ as above. Taking
$c=\g = \lfloor {\mathbb E}W-\s^2\rfloor$, so that the middle
term (almost) disappears, the expression can be rewritten as
\eq\label{f-ineq}
|{\mathbb E}[(W-\g)f(W)] - \l'{\mathbb E}[f(W+1)]|
  \leq \Bl\d + \sn b_i\Br\|\Delta f\|,
\en
where $\d$ and~$\l'$ are as above.

Fixing any set $B \subset \integ_+ + \g$, take $f = f_B$ as in~\Ref{f_B-def}.
It then follows from~\Ref{A.2} that
\eqa
\lefteqn{ |\pr(W\in B) - \TP(\ex W,\s^2)\{B\}| } \non \\
&=& |\ex\{(\bone_B(W) - \TP(\ex W,\s^2)\{B\})(I[W\ge \g] + I[W < \g])\}|
  \phantom{HH} \non\\
&\le& |\ex\{(\l' f_B(W+1) - (W-\g)f_B(W))\,I[W\ge \g]\}| + \pr(W < \g)\non\\
&=& |\ex\{\l' f_B(W+1) - (W-\g)f_B(W)\}| + \pr(W < \g),
\label{prob-diff}
\ena
this last from~\Ref{A.3}.
Hence \Ref{f-ineq} and~\Ref{prob-diff} show that, for any $B \subset\integ_+
+ \g$,
\eqa
\lefteqn{|\pr(W\in B) - \TP(\ex W,\s^2)\{B\}|}\non\\
&\le& \Bl\d + \sn b_i\Br \|\D f_B\| + \pr(W < \g) \non\\
&\le& (\l')^{-1}\Bl \d + \sn b_i\Br + \pr(W < \g).\label{A.6}
\ena
Now the largest value~$D$ of the differences $\{\TP(\ex W,\s^2)\{C\} -
\pr(W\in C)\}$, $C\subset\integ$, is attained at a set
$C_0 \subset \integ_++\g$,
and is thus bounded as in~\Ref{A.6}; the minimum is attained at~$\integ
\setminus C_0$ with the value~$-D$.  Hence
\[
|\pr(W\in C) - \TP(\ex W,\s^2)\{C\}|
          \le (\l')^{-1}\Bl \d + \sn b_i\Br + \pr(W < \g)
\]
for all $C\subset\integ$, and the lemma follows.
\ep

\bsk
If the random variables~$Y_i$ have finite variances,
both $\lambda'$ and~$\var W$ are typically of order~$O(n)$, so that
letting \hbox{$\bar{b}:= n^{-1}\sn b_i$}
 and applying Chebyshev's inequality to bound the final probability, we
find that then \hbox{$\|\law(W)-\TP(\ex W,\s^2)\|$} is of order
$O(n^{-1}+\bar{b})$.  Hence we are interested in choosing
$a_1, a_2, \ldots$ so that $b_1, b_2, \ldots$
are small. For independent $Y_1, Y_2, \ldots$,
 it is easy to convince oneself that the choice
\begin{equation}\label{a-choice}
a_i={\mathbb E}[Y_i\, W]-{\mathbb E}[Y_i]{\mathbb E}[W],
\end{equation}
is a good one, and this also emerges in our Markovian context.
 Notice that~\Ref{a-choice} implies that
$\s^2=\var W$.

Establishing~\Ref{i-bd} in the Markovian setting,
for~$a_i$ chosen as in~\Ref{a-choice}, is the core of the paper; it is accomplished
in Section~\ref{main}. For the estimates made in that analysis,
it is useful to introduce a coupling of $X$ with an independent copy
$X'=(X'_i)_{i=0}^\infty$. The relevant properties of the coupling
are given in the next section. From now on, we assume that
the conditional distributions $\law(Y_1 \giv X_1=k)$, $0\le k\le K$,
each have finite variance.

\section{The Markov chain coupling}\label{markov}
 \setcounter{equation}{0}
Let $X=(X_i)^\infty_{i=0}$ and $X'=(X'_i)^\infty_{i=0}$ be independent
copies
of an aperiodic, irreducible and stationary Markov chain with state
space $E=\{0,1,\ldots, K\}$. To understand their crucial role,
recall~\Ref{i-bd}, and note that
\begin{eqnarray}
{\mathbb E}[Y_if(W)]-{\mathbb E}[Y_i]{\mathbb E}[f(W)]
&=& {\mathbb E}[Y_if(W)]-{\mathbb E}[Y_if(W')]\nonumber\\
&=& {\mathbb E}[Y_i(f(W)-f(W'))].\label{mean-diff}
\end{eqnarray}
Here $W'=\sn Y'_i$, and $Y'_1,  \ldots, Y'_n$ are chosen
from the conditional distributions $(\law(Y_i \giv X'_i),\,1\le i\le n)$,
independently of each other and of $X$ and~$Y := (Y_1,\ldots,Y_n)$.
 Also, recall~\Ref{a-choice}, and note that then
\begin{equation}\label{a-equality}
a_i={\mathbb E}[Y_i (W-W')].
\end{equation}
Of course, \Ref{mean-diff} and~\Ref{a-equality} follow from the
independence of $(X,Y)$ and $(X',Y')$.

We refer to Lindvall~(2002, Part II.1) for proofs of the statements to
be made now; we shall be brief.

Let 0 be our reference state, and let $S=(S_m)^\infty_{m=0}$ and
$S'=(S'_m)^\infty_{m=0}$ be the points in increasing order of the sets
\[
\{k\in\integ_+\,;\, X_k=0\} \mbox{ and } \{k\in\integ_+\,;\, X'_k=0\},
\]
respectively. Then $S$ and $S'$ are stationary renewal
processes. Define $Z_0, Z_1,\ldots$, $Z'_0, Z'_1,\ldots$ by
\[
S_m=\sjom Z_j, \quad S'_m=\sjom Z'_j.
\]
Then all the $Z$ variables are independent, and the recurrence times $Z_1,
Z'_1, Z_2, Z'_2,\ldots$ are identically distributed, while the delays
$Z_0, Z'_0$ have the well-known distribution that renders $S$ and $S'$
stationary.

Now define $\tilde{S}=(\tilde{S}_m)^\infty_{m=0}$ to be the time points
at which both $S$ and~$S'$ have a renewal, i.e.
\[
\{k\in\integ_+\,;\, X_k=X'_k=0\}.
\]
Then $\tilde{S}$ is again a stationary renewal process, and we set
$\tilde{S}_m=\sjom \tilde{Z}_j$.

Let $X^\ast=(X^\ast_i)^\infty_{i=0}$ be an irreducible, finite state space
Markov chain with reference state 0, and let the associated
$(S^\ast_m)^\infty_{m=0}, (Z^\ast_j)^\infty_{j=0}$ have the obvious
meanings. For $j\geq 0$, write
\[
D_j=\min \{S^\ast_m-j;\, S^\ast_m \ge j\}.
\]
Due to the finiteness of the state space, it is easily proved that
there exists a $\rho > 1$ such that, as $m\to\infty$,
\begin{eqnarray}
&&  \max_k{\mathbb P}(D_j\geq m \giv X^\ast_j=k)=O(\rho ^{-m});
\label{ergod-1}\\
&&  {\mathbb P}(Z^\ast_0 \geq m)=O(\rho ^{-m}), \quad
 {\mathbb P}(Z^\ast_1 \geq m)=O(\rho ^{-m});\label{ergod-2}
\end{eqnarray}
c.f.~Lindvall~(2002, II.4, p.~30~ff.).
Of course, the maximum in~\Ref{ergod-1} does not depend on $j$.
When applied to
$((X_i, X'_i))^\infty_{i=0}$, the state space is $E \times E$; notice
that the aperiodicity of $X$ is needed to make $((X_i,
X'_i))^\infty_{i=0}$ irreducible.

For the rest of this section, drop the assumption that $X$ and $X'$
are stationary, but rather let $X_0=X'_0=0$, denoting the associated
probability by~$\pr^0$. We shall have much use
for an estimate of
\begin{equation}\label{8}
\b(n) :=
\left\|{\mathbb P}^0\left[\Bl X_n, 1+\sn Y_i \Br \in \cdot\right]
    -{\mathbb P}^0\left[\Bl X_n, \sn Y_i \Br \in \cdot\right] \right\|.
\end{equation}
It is natural to conjecture that
$\beta(n)=O(1/\sqrt{n})$, since that would be true if the sums $\sn Y_i$
formed a random walk independent of~$X$, under an aperiodicity assumption:
cf.~Lindvall~(2002, II.12 and II.14).

Let us say that the distribution of an integer-valued variable $V$ is
{\it strongly aperiodic\/} if
\eq\label{stap}
{\rm g.c.d.} \{k+i;\ {\mathbb P}(V=i)>0\}=1 \mbox{ for all } k.
\en
It is crucial to our argument to assume as smoothness condition that
\begin{eqnarray}\label{9}
& & \mbox{the distribution of } \sum^{S_1}_{i=1} Y_i
  \mbox { is strongly aperiodic,}
\end{eqnarray}
a condition that we are actually able to weaken later: see Theorem~\ref{Th-3}.
It then follows from~\Ref{9} that also
\begin{equation}\label{10}
\sum^{\tilde{S}_1}_{i=1} (Y_i-Y'_i) \mbox { is strongly aperiodic}.
\end{equation}
For the estimate of~\Ref{8}, notice that
\begin{eqnarray}
\lefteqn{\left\|{\mathbb P}^0\left[\Bl X_n,1+\sn Y_i \Br \in \cdot\right]
    -{\mathbb P}^0\left[\Bl X_n,\sn Y_i \Br \in \cdot\right] \right\|}\non\\
&=& \left\|{\mathbb P}^0\left[\Bl X_n,1+\sn Y_i \Br \in \cdot\right]
    -{\mathbb P}^0\left[\Bl X'_n,\sn Y'_i \Br \in \cdot\right] \right\|.
\label{11}
\end{eqnarray}
Now let
\[
\tau=\min \Blb k;\, 1+\sum^{\tilde{S}_k}_{i=1} Y_i
  =\sum^{\tilde{S}_k}_{i=1} Y'_i \Brb.
\]
We note that $\sum^{\tilde{S}_k}_{i=1} (Y_i-Y'_i)$, $k\ge0$,
is a random walk, with
step size distribution given by~\Ref{10}: it has expectation 0, finite
second moment, and is strongly aperiodic. For such a random walk,
Karamata's Tauberian theorem may be used to prove
that the probability that at least~$m$ steps are needed to hit
the state~$-1$ is of magnitude $O(1/\sqrt{m})$ (Breiman 1968, Theorem~10.25),
and hence
\begin{equation}\label{12}
{\mathbb P}^0(\tau\geq m)=O(1/\sqrt{m}).
\end{equation}
Now make a coupling as follows:
\[
X''_i=\left \{ \begin{array}{lll}
X'_i & \mbox{for} & i<\tilde{S}_\tau\\
X_i & \mbox{for} & i\geq \tilde{S}_\tau,
\end{array}\right .
\]
and define $Y''_i$, $i\geq 0$, accordingly. Recall~\Ref{11}. Standard
coupling arguments yield
\begin{eqnarray}
\lefteqn{\left\|{\mathbb P}^0\left[\Bl X_n,1+\sn Y_i\Br \in \cdot\right]
    -{\mathbb P}^0\left[\Bl X_n,\sn Y_i \Br \in \cdot\right] \right\|}\non\\
&=& \left\|{\mathbb P}^0\left[\Bl X_n,1+\sn Y_i \Br \in \cdot\right]
    -{\mathbb P}^0\left[\Bl X''_n,\sn Y''_i \Br \in \cdot\right] \right\|
          \non\\
&\le&  2{\mathbb P}^0(\tilde{S}_\tau >n).\label{13}
\end{eqnarray}
Let $\tilde{\mu}={\mathbb E}[\tilde{S}_1]$ and
$\alpha=1/(2\tilde{\mu})$. We get
\eqs
 {\mathbb P}^0(\tilde{S}_\tau>n)
&=& {\mathbb P}^0(\tilde{S}_\tau>n,\, \tau \geq \alpha n)+{\mathbb
  P}^0(\tilde{S}_\tau>n,\, \tau<\alpha n)\\
& \leq& {\mathbb P}^0(\tau \geq \alpha n)
  +{\mathbb P}^0(\tilde{S}_{\lfloor \alpha n \rfloor+1}>n).
\ens
But the latter probability is of order $O(1/n)$, due to Chebyshev's
inequality, and the former of order~$O(1/\sqrt n)$, by~\Ref{12}.
Hence \Ref{8}, \Ref{11} and~\Ref{13} imply that
\begin{equation}\label{beta(n)}
\beta (n)=O(1/\sqrt{n}) \mbox{  as  } n\rightarrow \infty.
\end{equation}

\section{Main theorem}\label{main}
 \setcounter{equation}{0}
We now turn to the approximation of~$\law(W)$, with~$W$ as defined in
the Markovian setting introduced in Section~\ref{intro}; the notation
is as in the previous section, and the
assumption that $X$ and~$X'$ are stationary is back in force.

In order to state the main lemma, we need some further terminology.
For each $1\leq i\leq n$, we define
\[
T^+_i:=\min\bigl\{n,\min\{\tilde{S}_k; \tilde{S}_k\geq i\}\bigr\}
\]
and
\[
T^-_i:=\left \{ \begin{array}{lll}
\max \{\tilde{S}_k; \tilde{S}_k\leq i\} & \mbox{if} & \tilde{S}_0\leq
i;\\
1 & \mbox{if} & \tilde{S}_0>i.
\end{array}\right .
\]
We then set
\eq\label{A-W-defs}
A_i=\sum^{T^+_i}_{j=T^-_i} Y_j,\quad
W^-_i =\sum^{T^-_i-1}_{j=1} Y_i,\ \mbox{ and }\
W^+_i  = \sum^n_{j=T^+_i+1} Y_j,
\en
with the understanding that $W^-_i=0$ if $T^-_i = 1$ and $W^+_i=0$
if $T^+_i = n$. We also define $A'_i$, ${W_i'}^{-}$  and ${W_i'}^{+}$
by replacing $Y_j$ by $Y'_j$.  For use in the argument to come, we
introduce independent copies $X\ul$ of the $X$-chain, $0\le l\le K$,
with $\law(X\ul) = \law(X \giv X_0=l)$.  By sampling the
corresponding $Y$-variables conditional on the realizations~$X\ul$,
we then construct the associated partial sum processes~$U\ul$ by setting
$U\ul_m := \sum_{s=1}^m Y\ul_s$.  Similarly, we define pairs of
processes $(\bX\ul,\bU\ul)$ in the same way, but based on
the {\it time-reversed\/} chain~$\bX$ starting with $\bX_0=l$
(Norris~1997, Theorem~1.9.1).  We use $\bbar(\cdot)$ to denote
the quantity in~\Ref{8} derived from the reversed chain, and note
that, under~\Ref{9}, the order estimate~\Ref{beta(n)} is true also
for~$\bbar$.
For any $m\ge1$ and $0\le l\le K$, we then write
\[
h_r(l,m) := \pr(U\ul_m \ge r+1),\ r\ge0;\qquad
h_r(l,m) := -\pr(U_m\ul \le r),\ r < 0,
\]
and specify $\bh_r(l,m)$ analogously, using the time-reversed processes~$\bU\ul$;
we then set $H(m) := \max\Blb \srz \|h_r(\cdot,m)\|,
\srz \|\bh_r(\cdot,m)\|\Brb$.

\begin{lemma}\label{main-este}
With the $a_i$ chosen as in~\Ref{a-choice}, the inequality~\Ref{i-bd}
is satisfied with
\eqs
b_i &:=& \tb(n/4)\Blb \half\ex\{|Y_i(A_i - A_i')|(|A_i|+|A_i'|)\}\vphantom{T_i^+A_I'}\right.\\
&&\left. \qquad\mbox{}
   + \ex\{|Y_i(A_i-A'_i)|(H(T_i^+-i) + H(i-T_i^-))\}\right.\\
&&\qquad\quad \left.\mbox{}
  + \ex\{|Y_i(A_i-A'_i)|\}\{\ex |A_i| + \ex(H(T_i^+-i) + H(i-T_i^-))\}\Brb\\
&&\quad\mbox{} + \g_i,
\ens
where, for $1\le i \le n/2$,
\eqs
   \g_i &=&  2 \ex\bigl\{|Y_i(A_i-A_i')|(I\tpnfg + \pr\tpnfg)\bigr\},\quad 1\le i \le n/2;\\
   \g_i &=&  2 \ex\bigl\{|Y_i(A_i-A_i')|(I\tmnfg + \pr\tmnfg)\bigr\},\quad  n/2 < i \le n,
\ens
and $\tb(n) = \max\{\b(n),\bbar(n)\} = O(n^{-1/2})$ under assumption~\Ref{9}.
\end{lemma}

\proof
The analysis of \Ref{i-bd} in our Markovian setting is
rather technical, and we divide it into three steps.
For the first step, we recall~\Ref{mean-diff}, giving
\begin{eqnarray}
\lefteqn{{\mathbb E}[Y_if(W)]-{\mathbb E}[Y_i]{\mathbb E}[f(W)]}\non\\
&=& {\mathbb E}[Y_i(f(W)-f(W'))]\non\\
&=& {\mathbb E}[Y_i(f(W^-_i+A_i+W^+_i) - f({W'}^{-}_i+A'_i+{W'}^{+}_i))]
   \non\\
&=& {\mathbb E}[Y_i(f(W^-_i+A_i+W^+_i)-f(W^-_i+A'_i+W^+_i))],
\label{15}
\end{eqnarray}
a careful proof of the last equality making use of a conditioning on
\[
\sigma\{T^-_i, T^+_i, \mbox{ and } X_j, X'_j, Y_j, Y'_j \mbox{ for }
T^-_i\leq j \leq T^+_i\},
\]
and of the symmetry of $X$ and $X'$.  Hence our aim is to bound
\eq\label{new-2}
   |{\mathbb E}\{Y_i(f(W^-_i+A_i+W^+_i)-f(W^-_i+A'_i+W^+_i)) - Y_i(A_i-A_i')\ex\D f(W)\}|\,.
\en

We now first consider indices~$i$ such that $1 \le i\le n/2$, and begin by observing 
that, by direct argument,  
\eqa
   &&{\mathbb E}\{[Y_i(f(W^-_i+A_i+W^+_i)-f(W^-_i+A'_i+W^+_i)) \non\\
   &&\qquad\qquad \qquad\mbox{} - Y_i(A_i-A_i')\ex\D f(W)]  I\tpnfg\} \non\\
   &&\qquad \Le 2\ex\{|Y_i(A_i-A_i')|I\tpnfg\}\,\|\D f\|. \label{new-1}
\ena
This brings part of the contribution to the quantity~$\g_i$ in the lemma, and allows
us to make the remaining argument assuming that $T_i^+-i \le n/4$. So let
\[
{\mathcal F}_i=\sigma\{T^-_i, T^+_i, \mbox{ and } X_j, X'_j, Y_j, Y'_j
\mbox{ for } 0\leq j\leq T^+_i\},
\]
and write
\eqa
  \lefteqn{{\mathbb E}\{Y_i(f(W^-_i+A_i+W^+_i)-f(W^-_i+A'_i+W^+_i))I\tpnfl\}} \label{16}\\
  &=&
   {\mathbb E}\{Y_iI\tpnfl\,{\mathbb E}[f(W^-_i+A_i+W^+_i)-f(W^-_i+A'_i+W^+_i)|{\mathcal F}_i]\}.
    \non
\ena
Now, for $r,m \in\integ$, define
\[
V(r,m) := I[0\le r \le m-1] - I[-1 \ge r\ge m],
\]
and observe that
\begin{eqnarray}
\lefteqn{f(W^-_i+A_i+W^+_i)-f(W^-_i+A'_i+W^+_i)}\non\\
&=& \srz \D f(W^-_i+A'_i+W^+_i + r)V(r,A_i-A_i')\non\\
&=& (A_i-A_i')\D f(W_i^- + W_i^+) \non\\
&&\quad\mbox{} +
  \srz [\D f(W^-_i+A'_i+W^+_i + r) - \D f(W_i^- + W_i^+)]
      V(r,A_i-A_i')\non\\
&=& (A_i-A_i')\D f(W_i^- + W_i^+) \non\\
&&\quad\mbox{} +
  \srz V(r,A_i-A_i')  \ssz \D^2f(W_i^- + W_i^+ + s) V(s,A_i'+r).
\label{17}
\end{eqnarray}
So, from \Ref{15}--\Ref{17}, 
we have isolated the term
\begin{equation}\label{18}
{\mathbb E}\{Y_i(A_i-A'_i)I\tpnfl\Delta f(W^-_i+W^+_i)\},
\end{equation}
from~\Ref{16},
together with an error involving the second differences in~\Ref{17}.
The remainder of the first step consists of bounding the
magnitude of this error.

To do so, note that the
second differences in~\Ref{17} are all of the form $\Delta^2f(\cdot\, +\,W^+_i)$,
where the ``$\cdot$''-part is measurable with respect to~${\mathcal F}_i$,
and~$W^+_i$ is the
contribution from the Markov chain starting from 0 at time~$T^+_i$. 
Furthermore, for any integer valued random variable~$Z$ and any bounded function~$h$,
$|\ex\D h(Z)| \le \|h\|\,\|\law(Z+1) - \law(Z)\|$.  Hence,
using~\Ref{8}, it follows that, for $T_i^+ - i \le n/4$,
entailing $n - T_i^+ \ge n/4$, we have
%
%
\begin{eqnarray}
|{\mathbb E}[\Delta^2f(\cdot +W^+_i)\giv{\mathcal F}_i]|
\,\leq\, \|\Delta f\|\beta(n/4),\hspace{20pt}\label{20}
\end{eqnarray}
where~$\b(n)$ is of magnitude
$\Ord{1/\sqrt{n}}$ under assumption~\Ref{9}, in view of~\Ref{beta(n)}.

Now observe that all the variables in~\Ref{17} except~$W^+_i$ are
${\mathcal F}_i$-measurable. What remains in order to use~\Ref{16}
is a careful count of the second difference terms in~\Ref{17}, of
which there are at most
$\frac{1}{2}|A_i-A'_i|(|A_i|+|A'_i|)$.
Using \Ref{15}--\Ref{18} and~\Ref{20}, we have found that
\eqa
    &&|{\mathbb E}\{Y_iI\tpnfl(f(W^-_i+A_i+W^+_i)-f(W^-_i+A'_i+W^+_i))]\}\non\\
  &&\qquad\mbox{} - 
       {\mathbb E}\{Y_i(A_i-A'_i)I\tpnfl\,\Delta f(W^-_i+W^+_i)\}|  \non\\
&&\quad \Le \half\b(n/4)\|\Delta f\|\,{\mathbb E}[|Y_i(A_i-A'_i)|(|A_i|+|A'_i|)]
\label{item-1}
\ena
where $\b(n)=O(1/\sqrt{n})$. This is responsible for the first
term in the expression for~$b_i$ in the statement of the lemma,
and completes the proof of the first step.

The next step is to work on $\ex[Y_i(A_i-A'_i)I\tpnfl\D f(W^-_i+W^+_i)]$.
 Although the random
variables $Y_i(A_i-A'_i)$, $W_i^-$ and $W_i^+$ are dependent, they
are conditionally independent given $T_i^-$ and~$T_i^+$, and then
$\law(W_i^+ \giv T_i^+=s) = \law(U\uo_{n-s})$ for $i\le s\le n$,
and $\law(W_i^- \giv T_i^- = s) = \law(\bU\uo_{s-1})$ for $1\le s\le i$.
This suggests writing
\eqa
\lefteqn{\ex[Y_i(A_i-A'_i)I\tpnfl\D f(W^-_i+W^+_i)]}\non\\
&=&
 \ex[Y_i(A_i-A'_i)I\tpnfl\D f(\bU\uo_{i-1} + U\uo_{n-i})] + \h_i\non\\
&=& \ex[Y_i(A_i-A'_i)I\tpnfl]\ex\{\D f(\bU\uo_{i-1} + U\uo_{n-i})\} + \h_i,
\label{2-split}
\ena
with~$\h_i$ to be bounded.

We start by writing
\eqs
   &&\ex[Y_i(A_i-A'_i)I\tpnfl\D f(W^-_i+W^+_i)] \\
  &&\qquad\Eq  \ex\{\ex[Y_i(A_i-A'_i)I\tpnfl\D f(W^-_i+W^+_i)\giv\gg_i] \},
\ens
with $\gg_i := \s(W_i^-,Y_i(A_i-A'_i),T_i^+)$.  Now $Y_i(A_i-A'_i)I\tpnfl$
is $\gg_i$-measurable, and
\eqs
\lefteqn{
\ex\{\D f(W_i^- + U\uo_{n-i}) - \D f(W_i^- + W_i^+) \giv \gg_i\}}\\
&=&  \ex\Blb \srz \D^2f(W_i^- + U\uo_{n-T_i^+} + r)\,
   V(r,U\uo_{n-i} - U\uo_{n-T_i^+}) \Giv \gg_i\Brb \\
&=&  \srz\ex\Blb  \D^2f(W_i^- + U\uo_{n-T_i^+} + r)\,
   h_r(X\uo_{n-T_i^+},T_i^+-i) \Giv \gg_i\Brb,
\ens
where the last line follows because, conditional on~$X\uo_{n-T_i^+}$,
$U\uo_{n-i} - U\uo_{n-T_i^+}$ is independent of $W_i^-$
and~$U\uo_{n-T_i^+}$. This in turn implies that, on $T_i^+-i \le n/4$,
\eqa
  \lefteqn{
  |\ex\{\D f(W_i^- + U\uo_{n-i}) - \D f(W_i^- + W_i^+) \giv \gg_i\}|}\non \\
&\le& \|\D f\| \srz\|h_r(\cdot,T_i^+-i)\|\,\non\\
&&\hspace{1.3cm}\times\, \ex \Blb \|\law(( X\uo_{n-T_i^+},U\uo_{n-T_i^+}+1))
    - \law(( X\uo_{n-T_i^+},U\uo_{n-T_i^+}))    \| \giv \gg_i \Brb\non\\
&\le& \|\D f\| H(T_i^+-i) \b(n/4), \label{new-3}
\ena
where the last line uses~\Ref{8}.  Thus it follows
that
\eqa
\lefteqn{
   |\ex\{Y_i(A_i-A'_i)I\tpnfl\D f(W_i^- + W_i^+)\}} \non\\
  &&\qquad\mbox{} - \ex\{Y_i(A_i-A'_i)I\tpnfl\D f(W_i^- + U\uo_{n-i})\}|\non\\
   &\le& \|\D f\|\b(n/4)\ex\{Y_i|A_i-A'_i|   H(T_i^+-i) \} =: \h_{i1}.
  \phantom{HHHHHHHHH}
\label{eta-i1}
\ena
An analogous argument, replacing $W_i^-$ by~$\bU\uo_{i-1}$, uses
the expression
\eqs
\lefteqn{
\ex\{\D f(\bU\uo_{i-1} + U\uo_{n-i}) - \D f(W_i^- + U\uo_{n-i})
    \giv \gg'_i\}}\\
&=&  \srz\ex\Blb  \D^2f(\bU\uo_{T_i^- -1} + U\uo_{n-i} + r)\,
   \bh_r(\bX\uo_{T_i^- -1},i-T_i^-) \Giv \gg'_i\Brb,
\ens
where $\gg'_i := \s(Y_i(A_i-A'_i),T_i^-,T_i^+)$, which we bound using~$\b(n-i)$
as a bound for $\|\law(U\uo_{n-i}+1)-\law(U\uo_{n-i})\|$, giving
\eqa
   \lefteqn{
     |\ex\{\D f(\bU\uo_{i-1} + U\uo_{n-i}) - \D f(W_i^- + U\uo_{n-i})
        \giv \gg'_i\}| } \non\\
  &\le& \|\D f\| \srz\|h_r(\cdot,i - T_i^-)\|\,
    \|\law(U\uo_{n-i}+1)-\law(U\uo_{n-i}) \| \non\\
  &\le& \|\D f\| H(i - T_i^-) \b(n/2). \label{new-4}
\ena 
This yields
\eqa
  &&|\ex\{Y_i(A_i-A'_i)I\tpnfl\D f(W_i^- + U\uo_{n-i})\} \non\\
  &&\qquad\quad\mbox{} - \ex\{Y_i(A_i-A'_i)I\tpnfl\D f(\bU\uo_{i-1} + U\uo_{n-i})\}| \non\\
  &&\quad\Le \|\D f\|\b(n/2)\ex\{Y_i|A_i-A'_i|   H(i-T_i^-) \} =: \h_{i2},
\phantom{HHHHHH}\label{eta-i2}
\ena
so that \Ref{2-split} holds with $\h_i = \h_{i1} + \h_{i2}$, accounting
for the second term in~$b_i$, and
completing the second step.

It now remains only to bound the difference between $\ex\{\D f(\bU\uo_{i-1}
+ U\uo_{n-i})\}$ and $\ex\{\D f(W)\}$.  This is accomplished much as
before, by writing $W = W_i^- + A_i + W_i^+$, and separating out the event
$T_i^+-i > n/4$.  This gives
\eqa
\lefteqn{|\ex\{(\D f(W) - \D f(W_i^- + W_i^+))I\tpnfl\}|}\non \\
&=& \Blm \ex\Blb \srz I\tpnfl \ex[\D^2f(W_i^- + W_i^+ + r)V(r,A_i)
   \giv T_i^-,T_i^+,A_i]\Brb\Brm  \non\\
&\le& \|\D f\|\b(n/4) \ex |A_i|,  \label{item-3}
\ena
where the last line is as for~\Ref{20}, and then
\eqa
  \lefteqn{|\ex\{\D f(W_i^- + W_i^+)I\tpnfl\} 
       - \ex\D f(\bU\uo_{i-1} + U\uo_{n-i}) \pr\tpnfl| }
     \non \\
  &&\Le \|\D f\| \b(n/4) \ex\{H(i - T_i^-) + H(T_i^+ - i)\},\phantom{HHHHHHHHHHHH}
\label{item-4}
\ena
this last as for \Ref{new-3} and~\Ref{new-4}. There is also the
inequality
\eqa
   &&|\ex\{\D f(W) I\tpnfg\} - \ex\{\D f(\bU\uo_{i-1} + U\uo_{n-i})\}\pr\tpnfg| \non\\
   &&\qquad  \Le 2\|\D f\|\, \pr\tpnfg,\label{new-5}
\ena
covering the contribution from $T_i^+-i > n/4$.
Multiplying the bounds in \Ref{item-3}, \Ref{item-4} and~\Ref{new-5}
by $\ex\{|Y_i(A_i - A'_i)|\}$ gives the third element of~$b_i$, together with 
the remaining contribution to~$\g_i$,
and the lemma is proved for $1\le i\le n/2$ .

For $n/2 < i \le n$, recall that $X$ and $X'$ are stationary. It
is well known
that then $(X_{n-j})^n_{j=0}$ and $(X'_{n-j})^n_{j=0}$ are also
stationary; these reversed processes inherit all the relevant properties
of $X$ and~$X'$. In carrying out the analysis above for the reversed
processes, we meet no obstacle, and hence the formula for the~$b_i$
holds also for $i > n/2$. This proves the lemma.
\ep

\msk
The bound in Lemma~\ref{main-este} can be combined with
Lemma~\ref{Stein-approx} to prove the total variation approximation
that we are aiming for, under appropriate conditions.  The expression
for~$b_i$ simplifies substantially, if we assume that
\eq\label{Y-assn}
\max\{\pr(Y_1 \ge r \giv X_1 = l), \pr(Y_1 \le -r \giv X_1=l)\}
  \Le \pr(Z \ge r)
\en
for all $r \ge 0$ and $0\le l\le K$,
for a positive integer valued random variable~$Z$ with $\ex Z^3 < \infty$.
If this is the case, then
\eqs
H(m) \le 2m\ex Z;&&
   \ex(|A_i|\giv X,X') \le 2(T_i^+ - T_i^- + 1)\ex Z;\\
\ex(A_i^2\giv X,X') &\le&  2(T_i^+ - T_i^- + 1)^2\ex Z^2;\\
\ex(|Y_iA_i| \giv X,X') &\le& 2(T_i^+ - T_i^- + 1)\ex Z^2
\ens
and
\eqs
\ex(|Y_i|A_i^2 \giv X,X') &\le&  2(T_i^+ - T_i^- + 1)^2\ex Z^3  .
\ens
{}From these bounds, together with the
fact that $A_i$ and~$A'_i$ are independent conditional on $X,X'$,
it follows that
\eqa
b_i &\le& \tb(n/4)\{4\ex Z^3\ex\t_i^2 + 8\ex Z\ex Z^2\ex\t_i^2
         + 4\ex Z^2\ex\t_i(2\ex Z\ex \t_i + 2\ex Z\ex \t_i)\}\non\\
  &&\qquad + 8(C\vee \overline C)(n \ex Z^2 \r^{-n/4} + \ex\t_i\ex Z^2 \r^{-n/4}) \non\\
&\le& 28\tb(n/4) \ex Z^3\ex\t_i^2 + 16(C\vee \overline C)n\ex Z^2 \r^{-n/4}, \label{b_i-simple}
\ena
where $\t_i := T_i^+ - T_i^- + 1$, and $C,\overline C$ are the constants implied in~\Ref{ergod-1}
for the $X$--process and its time reversal.  Note that, since $X$ and~$X'$
are in equilibrium, both chains can be taken to run for all positive and
negative times, so that then $\ex\t_i^2 \le \ex\t^2$, where~$\t$ is the
length of that interval between successive times at which both $X$
and~$X'$ are in the state~$0$ which contains the time point~$0$.
$\ex\t_i^2$ is in general smaller than $\ex \t^2$,
because $T_i^-$ and~$T_i^+$
are restricted to lie between~$1$ and~$n$. Then the bound~\Ref{b_i-simple},
combined with
Lemma~\ref{Stein-approx}, leads to the following theorem.

\begin{theorem}\label{Th-1}
Under assumptions \Ref{9} and~\Ref{Y-assn}, and with stationary~$X$,
it follows that
\[
\|\law(W) - \TP(\ex W,\var W)\|
   \le 4\Bl 1 + 14n\f(n)\ex\t^2\ex Z^3\Br/\var W,
\]
where $\f(n) := \tb(n/4) + (C\vee \overline C)n\r^{-n/4}$.
\end{theorem}

\nin Note that
\[
\var W = \sn \ex [\var(Y_i \giv X_i)] + \var\Bl\sn\ex(Y_i\giv X_i)\Br,
\]
so that the bound in Theorem~\ref{Th-1} is of order
$\Ord{n^{-1} + \f(n)} = O(n^{-1/2})$ under these assumptions,
unless $\law(Y_1)$ is degenerate, in which case~$W$ is a.s.\
constant. Note also that replacing each~$Y_i$ by $Y_i-c$, for any
$c\in\integ$, results only in a translation, and does not change
$\|\law(W) - \TP(\ex W,\var W)\|$, and this can be exploited if
necessary when choosing the random variable~$Z$ in~\Ref{Y-assn}.

\msk
The assumption that~$X$ be stationary is not critical.

\begin{theorem}\label{Th-2}
Suppose that the assumptions of Theorem~\ref{Th-1} hold, except that
the initial distribution
$\law(X_0)$ is not the stationary distribution.
Then it is still the case that
$\|\law(W) - \TP(\ex W,\var W)\| = O(n^{-1/2})$.
\end{theorem}

\proof
Let~$X'$ be in equilibrium and independent of~$X$, and use it as
in Section~\ref{markov} to construct an equilibrium process~$X''$
which is identical with~$X$ after the time~$T_1^+$ at which $X$
and~$X'$ first coincide in the state~$0$.  Then Theorem~\ref{Th-1}
can be applied to~$W''$, constructed from~$X''$, and also
\[
W = A_1 + W_1^+ \ \mbox{ and }\ W'' = A_1'' + W_1^+,
\]
with $A_1$ and~$A_1''$ defined as before.  Let $g:\integ\to\re$ be
any bounded function, and observe that
\eqa
\lefteqn{|\ex g(W) - \ex g(W'')|
              = |\ex\{g(A_1+W_1^+) - g(A_1'' + W_1^+)\}|}\non\\
&\le& \Blm\ex\Blb \ex\Bl I[A_1 > A_1''] \sum_{j=1}^{A_1-A_1''}
   \D g(W_1^+ + A_1'' + j - 1)\giv T_1^+,A_1,A_1''\Br
     \right.\right.\hspace{40pt}\non\\
&&\ \left.\left.\mbox{}
     - \ex\Bl I[A_1 < A_1''] \sum_{j=1}^{A_1''-A_1}
         \D g(W_1^+ + A_1 + j - 1)\giv T_1^+,A_1,A_1''\Br\Brb \Brm.
\label{2.1}
\ena
Now, arguing as before, on $T_1^+ \le n/2$, we have
\[
|\ex\{\D g(W_1^++A''_1+j) \giv T_1^+,A_1,A''_1\}| \le \|g\|\,
  \|\law(W_1^+ +1) - \law(W_1^+)\| \le \|g\|\b(n/2),
\]
with $\b(n) = O(1/\sqrt n)$, implying from~\Ref{2.1} that
\eqa
|\ex g(W) - \ex g(W'')|  &\le& \{2\pr[T_1^+ > n/2] + \ex|A_1-A_1''|\b(n/2)\}\|g\| \non\\
   &\le& 4\ex T_1^+\ex Z\f(n)\|g\|. \label{New-star}
\ena
Although the distribution of~$T_1^+$ is not the same as if both
$X$ and~$X'$ were at equilibrium, it has moments which are
uniformly bounded for all initial distributions~$\n$, in view of
\Ref{ergod-1} and~\Ref{ergod-2},
and hence, from~\Ref{New-star} and because $\f(n) = O(1/\sqrt n)$,
it follows that $\|\law(W)-\law(W'')\| = O(n^{-1/2})$.

On the other hand,
\[
|\ex W - \ex W''| \Le \ex|A_1-A_1''| \Le 2\ex T_1^+\ex Z,
\]
and also
\[
|\var W - \var W''| \Le \var(W-W'') + 2\sqrt{\var W\,  \var(W-W'')},
\]
with
\[
  \var (W-W'') \Le \ex\{|A_1-A_1''|^2\} \Le 4 \ex \{(T_1^+)^2\}\ex \{Z^2\} \ =:\ 4D^2,
\]
giving
\[
  |\var W - \var W''| \Le 8 D \max\{\sqrt{\var W},D\}.
\]
Hence, from Lemma~\ref{TP-diff}, it follows that
\[
\|\TP(\ex W,\var W) - \TP(\ex W'',\var W'')\| = O(n^{-1/2})
\]
also, completing the proof.
\ep

\bsk
Assumption~\Ref{9}, that the distribution~$Q:=\law\Bl\sum_{i=1}^{S_1}
Y_i \Giv X_0=0 \Br$ be strongly aperiodic, can actually be relaxed;
it is enough to assume that~$Q$ is aperiodic.

\begin{theorem}\label{Th-3}
Suppose that the assumptions of Theorem~\ref{Th-1} hold, except that
assumption~\Ref{9} is weakened to assuming that~$Q$ is aperiodic.
Then it is still the case that
$\|\law(W) - \TP(\ex W,\var W)\| = O(n^{-1/2})$.
\end{theorem}

\proof
Define a new Markov chain~$\hX$ by splitting the state~$0$ in~$X$ into
two states, $0$ and~$-1$.  For each~$j$, set
\[
\hX_j = \begin{cases}
  X_j &\mbox{ if}\ X_j \ge 1;\cr
  -R_j &\mbox{ if}\ X_j=0,\cr
\end{cases}
\]
where $(R_j,\,j\ge0)$ are independent Bernoulli $\Be(1/2)$ random
variables; then set $\hY_j = Y_j$, $j\ge0$, and define $\hW = \sjn \hY_j$.
Clearly, $W=\hW$ a.s., so that we can use the construction based on the
chain~$\hX$ to investigate~$\law(W)$.  However, choosing~$0$ as reference
state also for~$\hX$, we have
\[
\hQ := \law\Bl \sum_{j=1}^{\hS_1} \hY_j \Giv \hX_0=0 \Br
  = \law\Bl \sum_{m=1}^M V_m \Br,
\]
where $V_1,V_2,\ldots$ are independent and identically distributed
with distribution~$Q$, and~$M$ is independent of the~$V_j$'s, and has
the geometric distribution~$\Ge(1/2)$. Since~$Q$ is aperiodic, it
follows that~$\hQ$ assigns positive probability to all large enough
integer values, and is thus strongly aperiodic.  Hence Theorems
\ref{Th-1} and~\ref{Th-2} can be applied to~$W$, because of its
construction as~$\hW$ by way of~$\hX$ and~$\hY$.
\ep

\msk
\Ceka\ \&~Mikalauskas~(1999) have also studied total variation
approximation in this context, in the degenerate case in which
$Y_1 = h(k)$ a.s.\ on $\{X_1=k\}$, $0\le k\le K$.  They use
characteristic function arguments, based on earlier work of
Sira\v zdinov \&~Formanov~(1979), and their approximations are
in terms of signed measures, rather than translated Poisson
distributions.  In their Theorem~2.2, they give one approximation
with error of order $O(n^{-1/2})$, and another, more
complicated approximation with error of order $o(n^{-1/2})$.
However, their formulation is probabilistically opaque, and
their proofs give no indication as to the magnitude of the
implied constants in the error bounds, or as to their dependence
on the parameters of the problem.  In fact, their `smoothness'
condition~(2.8) requires that the Markov chain~$X$ has a certain
structure, irrespective of the values of~$h$, which is unnatural.
For example, the $X$-chain with $K=2$ which has transition matrix
\eq\label{New-2}
\Bl
\begin{matrix}
\tfrac9{10}&\tfrac1{10}&0\\
         0&0&1\\
    1&0&0
\end{matrix}
\Br
\en
fails to satisfy their condition, although, for many score
functions~$h$, \Ref{1.3} is still true; for instance, our Theorem~\ref{Th-3}
applies to prove~\Ref{1.3} if $h(0)=3$ and $h(1)=h(2)=1$. However,
$Q$ is not aperiodic when $h(0)=3$, $h(1)=1$ and $h(2)=2$, and, without
this smoothness condition being satisfied, Theorem~\ref{Th-3} cannot
be applied.  This is in fact just as well, since the equilibrium
distribution of~$W$ then assigns probability much greater than
$\frac23$ to the set $3\integ \cup \{3\integ + 1\}$, whereas the
probability assigned to this set by the translated Poisson distribution
with the corresponding mean and variance approaches $\frac23$ as
$n\to\infty$.

\adb{In fact, if~$Q$ is periodic, it is rather the exception than
the rule that $\law(W)$ and $\TP(\ex W,\var W)$ should be close in total
variation.  To see this, let~$Q$ have period~$d$.  Fix
any $k\in E$, and take any~$i\in\ZZ_+$ and any realization of
the process such that $X_0(\o)=0$ and $X_i(\o)=k$; let
$R_{ki}(\o) := \sum_{l=1}^i Y_l(\o)$  modulo~$d$. Then it is immediate
that~$R_{ki}(\o) = r_k$ is a constant depending only on~$k$, since,
continuing two such realizations along the same $X$-path
and with the same~$Y$ values until the process next hits~$0$,
the two $Y$-sums then have to have the same remainder~$0$ modulo~$d$.
The same considerations show that $\law(Y_i \giv X_i=k)$ is
concentrated on a set $d\ZZ + \r_k$ for some $\r_k \in \{0,1,\ldots,d-1\}$,
and that the transition matrix~$P=(p_{kj})$ of the $X$-chain satisfies the
condition
\eq\label{P-cond}
r_k + \r_j \equiv r_j  \mod d \quad \mbox{whenever}\quad p_{kj} > 0.
\en
Moreover, for the same $r$- and $\r$-values,
{\it any\/} choice of~$P$ consistent with~\Ref{P-cond} yields
a distribution~$Q$ with period~$d$.}

\adb{Now the distribution $\TP(\m,\s)$ assigns probability approaching
$1/d$ as $\s\to\infty$ to any set of the form $d\ZZ + r$, $r\in
\{0,1,\ldots,d-1\}$. On the other hand, using $\pr^\l$ to denote
probabilities computed with~$\l$ as the distribution of~$X_0$, we have
\eqs
\pr^\l[W \equiv r \!\! \mod d] &=& \sie \l_i \pr[W\equiv r \!\!\mod d \giv X_0 = i]\\
  &=& \sie \l_i \pr[X_n \in E_{r-r_i} \giv X_0 = i],
\ens
where $E_r := \{k\in E: r_k = r\}$ and
differences in the indices are evaluated modulo~$d$.
This, as $n\to\infty$, approaches the value
\[
\sie \l_i \p(E_{r-r_i}) = \sodmi \l(E_s)\p(E_{r-s}),
\]
where~$\p$ is the stationary distribution of the $X$-chain.
Hence $\law^\l(W)$ becomes far from any
translated Poisson distribution as $n\to\infty$ unless
\eq\label{lambda-cond}
\sodmi \l(E_s)\p(E_{r-s}) = 1/d \ \mbox{for all}\ 0\le r\le d-1.
\en
}

\adb{
It is immediate that~\Ref{lambda-cond} cannot hold for {\it all\/}
choices of~$\l$ unless
\eq\label{pi-cond}
\p(E_r)=1/d \  \mbox{for each}\  r\in\{0,1,\ldots,d-1\}.
\en
What is more, it cannot hold in the stationary case, when $\l=\p$,
unless~\Ref{pi-cond} holds.  This follows from
multiplying both sides of~\Ref{lambda-cond} (with $\l=\p$) by~$t_j^r$ and
adding over~$r$, where~$t_j$, $0\le j\le d-1$, are the
complex $d$-th roots of unity, with $t_0 := 1$.
Writing $\p(t) := \sodmi \p(E_s) t^s$, this implies that $\{\p(t_j)\}^2 = 0$
for $1\le j\le d-1$,
and hence that the polynomial~$\p(t)$ is proportional to the polynomial
$\sodmi t^s$, which implies~\Ref{pi-cond}.
Indeed, the (circulant) matrix~$\Pi$ with elements $\Pi_{rs} =\p(E_{r-s})$
has~$d$ distinct eigenvectors corresponding to the eigenvalues~$\p(t_j)$,
so that
if $\p(t_j) \neq 0$ for all~$j$, then~\Ref{lambda-cond} has $\l(E_s)=1/d$
for all~$s$ as its only solution.}

\adb{
But condition~\Ref{P-cond}
depends only on the communication structure of~$P$, and not on
the exact values of its positive elements, whereas for~\Ref{pi-cond}
to be true
needs careful choice of the values of these elements.  Hence,
for most choices of~$P$ leading to a periodic~$Q$, meaning those in which
$\p(E_r) = 1/d$ for all~$r$ is {\it not\/} true, $\law^\l(W)$
and $\TP(\ex^\l W,\var^\l W)$ are not asymptotically close for $\l=\p$,
or if~$\l$ is concentrated on a single point, or indeed, if $\p(t_j)
\neq 0$ for all~$j$, for any~$\l$ not satisfying $\l(E_s)=1/d$
for all~$s$.  In consequence, for most choices of~$P$ leading to a
periodic~$Q$, the conclusions
of Theorems \ref{Th-1} and~\ref{Th-2} are very far from true.}

\adb{In the example~\Ref{New-2} above, $Q$ has period~$3$ when $h(0)=3$,
$h(1)=1$ \hbox{and~$h(2)=2$.}  Clearly, we have $\r_0=0$, $\r_1=1$
and~$\r_2=2$; we then also have $r_0=r_2=0$ and $r_1=1$, so that $E_0 = \{0,2\}$,
$E_1=\{1\}$ and $E_2=\emptyset$.  It is
easy to check that~\Ref{P-cond} is satisfied, and that it
would still be satisfied if $p_{21}$ were also positive. The
matrices~$P$ consistent with condition~\Ref{P-cond}
for these values of the $\r_k$ and~$r_k$ thus take the form
\[
\Bl
\begin{matrix}
1-\a&\a&0\\
0&0&1\\
\b&1-\b&0
\end{matrix}
\Br
\]
for $0 \le \a,\b \le 1$, so that $\p = (\b,\a,\a)/\{\b+2\a\}$;
in~\Ref{New-2}, $\a = 1/10$ and $\b=1$.
However, since~$\p(E_2)$ is necessarily zero, condition~\Ref{pi-cond}
is never satisfied. Furthermore, $\l(E_2)$ must also be zero, and
$\p(t_j)=0$ can only occur for~$t_j$ a complex
cube root of unity if $\b=\a$. Thus, in this example, the conclusions
of Theorems \ref{Th-1} and~\ref{Th-2} are never true; furthermore,
if $\a\neq\b$, translated Poisson approximation cannot be good
for {\it any\/} initial distribution~$\l$.}

\adb{As a second example, take $K=3$ and $P$ of the form
\[
\Bl
\begin{matrix}
1-\a&\a&0&0\\
0&0&1&0\\
0&0&1-\b&\b\\
1&0&0&0
\end{matrix}
\Br
\]
for $0 \le \a,\b \le 1$, so that $\p = (\b,\a\b,\a,\a\b)/\{\a+\b+2\a\b\}$.
This matrix satisfies~\Ref{P-cond} for $Y$-distributions satisfying
$\r_0 = \r_2 = 0$ and $\r_1 = \r_3 = 1$ with $d=2$, and then $r_0 = r_3 = 0$
and $r_1=r_2=1$, so that $E_0 = \{0,3\}$ and $E_1 = \{1,2\}$.  Hence
$\p(E_0)=\p(E_1)=1/2$ only if $\a=\b$, and, if $\a\neq\b$, $\p(-1)\neq0$.
Thus, if $\a\neq\b$, the conclusions of Theorems \ref{Th-1} and~\ref{Th-2}
are far from true, and indeed translated Poisson approximation cannot possibly
be good for any initial distribution~$\l$ which does not give equal weight
to $E_0$ and~$E_1$.}

The assumption that $X$ has finite state space~$E$ greatly
simplifies our arguments, because uniform bounds on hitting and coupling
times, such as those given in \Ref{ergod-1} and~\Ref{ergod-2}, are
immediate.  Results similar to ours can be expected to hold also
for countably infinite~$E$, provided that the chain~$X$ is such
that uniform bounds analogous to \Ref{ergod-1} and~\Ref{ergod-2} are
valid, and if the distributions of the~$Y_i$ are such that,
for instance, \Ref{Y-assn} also holds. However, a full analysis
of the case in which~$E$ is countably infinite would be
a substantial undertaking.

\section{Acknowledgements}
The authors wish to thank the referees for many helpful suggestions.
ADB is grateful for support from the Institute for Mathematical
Sciences of the
National University of Singapore, from Monash University School
of Mathematics and
Statistics, and from Schweizer Nationalfonds Projekt Nr~20--67909.02.

\end{document}